\documentclass[11pt]{article}
\usepackage[T2A]{fontenc}
\usepackage{amsmath}
\usepackage{amssymb}
\usepackage{arydshln}
\usepackage[colorlinks=true]{hyperref}
\def\a{\mathbf{a}}
\def\cb{\mathbf{c}}
\def\d{\mathbf{d}}
\def\b{\mathbf{b}}
\def\aaaa{a}
\def\bbbb{b}

\def\cccc{c}
\def\C{\mathbb{C}}

\def\G{\mathcal{I}}
\def\L{\mathcal{L}}
\def\GI{\mathfrak{G}}
\def\GO{\mathcal{G}}

\def\f{\mathbf{f}}
\def\m{\mathbf{m}}
\def\N{\mathbb{N}}
\def\R{\mathbb{R}}
\def\etta{\boldsymbol{\eta}}
\def\lambdda{\boldsymbol{\lambda}}
\def\zetta{\boldsymbol{\zeta}}

\def\ggamma{\boldsymbol{\gamma}}
\newtheorem{theorem}{\hspace*{\parindent}Theorem}
\newtheorem{lemma}{\hspace*{\parindent}Lemma}
\newtheorem{corollary}{\hspace*{\parindent}Corollary}

\DeclareMathOperator*{\res}{\mathrm{res}}
\title{On Meijer's $G$ function $G^{m,n}_{p,p}$ for $m+n=p$}
\author{D.B.\:Karp$^{\rm a,b}$\footnote{Corresponding author. E-mail: D. Karp -- \emph{dimkrp@gmail.com},
E.\:Prilepkina --  \emph{pril-elena@yandex.ru}}~~and
E.G.\:Prilepkina$^{\rm b,c}$
\\[10pt]
\\
\small{\textit{$\phantom{1}^a$Holon Institute of Technology, Holon, Israel}}
\\
\small{\textit{$\phantom{1}^b$Far Eastern Federal University, Ajax Bay~10, 690922 Vladivostok, Russia}}
\\
\small{\textit{$\phantom{1}^c$Institute of Applied Mathematics,
FEBRAS, 7 Radio Street, Vladivostok,  690041, Russia}}}
\date{}
\begin{document}
\maketitle


\begin{abstract}
The paper is devoted to the piece-wise analytic case of Meijer's $G$ function $G^{m,n}_{p,p}$.  While the problem of its analytic continuation was solved in principle by Meijer and Braaksma we show that in the ''balanced'' case $m+n=p$ the formulas take a particularly simple form.  We derive explicit expressions for the values of these analytic continuations on the banks of the branch cuts.  It is further demonstrated that particular cases of this type of $G$ function having integer parameter differences satisfy identities similar to the Miller-Paris transformations of the generalized hypergeometric function. Finally, we give a presumably new integral evaluation involving $G^{m,n}_{p,p}$  function with $m=n$ and apply it for summing a series involving digamma function and related to the power series coefficients of the product of two generalized hypergeometric functions with shifted parameters.       
\end{abstract}

\bigskip


\bigskip


\section{Introduction and preliminaries}

Let us start by setting the stage. Throughout the paper we will use the shorthand notation 
\begin{equation}\label{eq:notation}
\begin{split}
\Gamma(\a)&=\prod\nolimits_{i=1}^{p}\Gamma(a_i),~~
\sin(\a)=\prod\nolimits_{i=1}^{p}\sin(a_i)
\\
\a+\beta&=(a_1+\beta,\ldots,a_p+\beta),~~\a_{[k]}=(a_1,\ldots,a_{k-1},a_{k+1},\ldots,a_p),
\end{split}
\end{equation}
for any vector $\a=(a_1,\ldots,a_p)\in\C^{p}$ (the set of complex  $p$-tuples), a scalar $\beta$ and with $\Gamma(\cdot)$ standing for  Euler's gamma function.  The main character of this work is Meijer's $G$ function defined as follows.  Suppose $0\leq{m}\leq{q}$, $0\leq{n}\leq{p}$ are integers and
$\a=(\a_1,\a_2)$ with $\a_1\in\C^m$, $\a_2\in\C^{q-m}$,  $\b=(\b_1,\b_2)$, $\b_1\in\C^n$, $\b_2\in\C^{p-n}$ are complex vectors satisfying $a_i-b_j\notin\N$ for all $i=1,\ldots,m$ and $j=1,\ldots,n$.  Meijer's $G$-function is defined  by the Mellin-Barnes integral  
\begin{equation}\label{eq:G-defined}
G^{m,n}_{p,q}\!\left(\!z~\vline\begin{array}{l}\b\\\a\end{array}\!\!\right)\!\!:=
\\
\frac{1}{2\pi{i}}
\int\limits_{\mathcal{L}}\!\!\frac{\Gamma(\a_1\!+\!s)\Gamma(1-\b_1\!-\!s)}
{\Gamma(\b_{2}\!+\!s)\Gamma(1-\a_{2}\!-\!s)}z^{-s}ds,
\end{equation}
where the contour $\L$ is a simple loop that starts and ends at infinity and separates the poles of $s\to\Gamma(\a_1+\!s)$ leaving them on the left from those of $s\to\Gamma(1-\b_1-\!s)$ leaving them on the right. Details regarding the choice of the contour $\L$ and the convergence of the above integral can be found, for instance, in \cite[section~1.1]{KilSaig}, \cite[Section~8.2]{PBM3}, \cite[section~16.17]{NIST} and \cite[Appendix]{KPSIGMA2016}. 
The power function $z^{-s}$ is generally defined on the Riemann surface of the logarithm, so that
$$
z^{-s}=\exp(-s\{\log|z|+i\arg(z)\})
$$
and $\arg(z)$ is allowed to take any real value.  In this paper, however, we will restrict our attention to the principal value of the argument $-\pi<\arg(z)\le\pi$. $G$ function is a natural generalization of the generalized hypergeometric function and is indispensable when constructing the fundamental solutions of the generalized hypergeometric differential equation  \cite{KPSIGMA2016,Norlund,Scheidegger}. It is also important in probability and statistics \cite{Consul,ST}, random matrix theory and in adjacent field of multiple orthogonal polynomials \cite{BGS,KuijMolag}.

Denoting the integrand in \eqref{eq:G-defined} by 
$$
\G(s)=\frac{\Gamma(\a_1\!+\!s)\Gamma(1-\b_1\!-\!s)}
{\Gamma(\b_{2}\!+\!s)\Gamma(1-\a_{2}\!-\!s)},
$$
we have the following theorem by combining the known existence results \cite[Theorems~1.1,1.2]{KilSaig} or \cite[Theorem~A.1]{KPSIGMA2016}:
\begin{theorem}
Suppose $m+n\le p$. Then the function $G^{m,n}_{p,p}$ is piece-wise analytic with 
\begin{equation}\label{eq:sumresleft}
G^{m,n}_{p,p}\left(z\left|\begin{array}{l} \b\\
\a\end{array}\right.\right)=\sum\limits_{j=1}^m\sum\limits_{l=0}^\infty\res_{s=-a_j-l}\G(s)z^{-s}
\end{equation}
for $|z|<1$ and 
\begin{equation}\label{eq:sumresright}
G^{m,n}_{p,p}\left(z\left|\begin{array}{l} \b\\
\a\end{array}\right.\right)=-\sum\limits_{i=1}^n\sum\limits_{k=0}^\infty\res_{s=1-b_i+k}\G(s)z^{-s}
\end{equation}
for $|z|>1$.
\end{theorem}

The functions obtained by summing the  series \eqref{eq:sumresleft} and \eqref{eq:sumresright} are, generally speaking, not analytic continuations of each other.   On the other hand, we can continue the function defined by  \eqref{eq:sumresleft} inside the unit disk to the domain $|z|>1$ (the analytic continuation NOT being given by \eqref{eq:G-defined} or \eqref{eq:sumresright}) and the function defined by  \eqref{eq:sumresright} outside the unit disk to the domain $|z|<1$ (again the analytic continuation is NOT given by \eqref{eq:G-defined} or \eqref{eq:sumresleft}).  Hence, it makes sense to introduce a separate notation for the function defined by \eqref{eq:G-defined} and \eqref{eq:sumresleft} for $|z|<1$ and analytically continued to $|z|>1$.  We will call this function ''internal'' and denote it by  
$\GI^{m,n}_{p,p}(z)$. Similarly, the function defined by \eqref{eq:G-defined} and \eqref{eq:sumresright} for $|z|>1$ and analytically continued to $|z|<1$ will be call ''external'' and denoted by $\GO^{m,n}_{p,p}(z)$. In view of this definition, the standard $G$ function $G^{m,n}_{p,p}(z)$ is given by:
$$
G^{m,n}_{p,p}\left(z\left|\begin{array}{l}\b\\\a\end{array}\right.\right)
=\GI^{m,n}_{p,p}\left(z\left|\begin{array}{l}\b\\\a\end{array}\right.\right)~\text{if}~|z|<1;
~~G^{m,n}_{p,p}\left(z\left|\begin{array}{l}\b\\\a\end{array}\right.\right)
=\GO^{m,n}_{p,p}\left(z\left|\begin{array}{l}\b\\\a\end{array}\right.\right)~\text{if}~|z|>1.
$$
The the well-known reflection property of $G$-function takes the form:
\begin{equation}\label{eq:reflection}
\GI^{m,n}_{p,q}\!\left(z~\vline\begin{array}{l}\b\\\a\end{array}\!\right)
\!=\!\GO^{n,m}_{q,p}\!\left(\frac{1}{z}~\vline\begin{array}{l}1-\a\\1-\b\end{array}\!\right)~\text{if}~|z|<1;
~\GO^{m,n}_{p,q}\!\left(z~\vline\begin{array}{l}\b\\\a\end{array}\!\right)
\!=\!\GI^{n,m}_{q,p}\!\left(\frac{1}{z}~\vline\begin{array}{l}1-\a\\1-\b\end{array}\!\right)~\text{if}~|z|>1.
\end{equation}
If all  poles of the integrand $\G(s)$ are simple, then the series  in \eqref{eq:sumresleft} and \eqref{eq:sumresright}, and hence each of the functions  $\GI^{m,n}_{p,p}(z)$, $\GO^{m,n}_{p,p}(z)$ can be written as a finite sum of hypergeometric functions. This expansion was derived by Cornelis Meijer himself, see \cite[(7)]{Meijer} and further details in \cite[section~4.6.2]{Slater}.  Then one can use the well known analytic continuation formula for the generalized hypergeometric function to construct the analytic continuations of $\GI^{m,n}_{p,p}(z)$ and $\GO^{m,n}_{p,p}(z)$. The  problem of analytic continuation was tackled by Meijer \cite[Theorem~22,p.1169]{Meijer} and Braaksma \cite[Theoreom~10,p.320]{Braaksma}.  A brief summary of their results is as follows.  Combination of \cite[(214), (31), (15), (88)]{Meijer} yields  ($m+n-p=\varkappa$ for brevity)
\begin{multline}\label{eq:MeijerContin}
G^{m,n}_{p,p}\!\left(\!z~\vline\!\begin{array}{l}\b\\\a\end{array}\!\!\right)
\!=\!-\frac{(-\pi)^{\varkappa}}{z}\sum\limits_{k=1}^{p}\frac{\Gamma(b_k-\b_{[k]})}{\Gamma(b_k-\a)}z^{b_k}
{_pF_{p-1}}\!\left(\begin{array}{l}1-\a+b_k\\1-\b_{[k]}+b_k\end{array}\vline\:\frac{(-1)^{\varkappa}}{z}\right)
\\
\times\sum\limits_{j=1}^{m}\frac{\sin(\pi(\b_2-a_j))e^{i\pi(\varkappa+2l-1)(a_j-b_k)}}{\sin(\pi(\a_{1{[j]}}-a_{j}))\sin((a_j-b_k)\pi)}
\end{multline}
for $(\varkappa+2l-2)\pi<\arg(z)<(\varkappa+2l)\pi$ and under additional restrictions that the components of $\b$ and $\a_1$ are distinct modulo integers.  For the same range of $\arg(z)$ Braaksma's formula \cite[(11.38)]{Braaksma} after correcting a misprint and a small error reads as follows 
\begin{equation*}
G^{m,n}_{p,p}\!\left(\!z~\vline\!\begin{array}{l}\b\\\a\end{array}\!\!\right)
\!=\!-S(z)+\sum\limits_{j=1}^{m}\frac{\sin(\pi(\b_2-a_j))e^{i\pi(\varkappa+2l-1)a_j}}{\pi\sin(\pi(\a_{1{[j]}}-a_{j}))}
G^{1,p}_{p,p}\!\left(\!ze^{-i\pi(\varkappa+2l-1)}~\vline\!\begin{array}{c}\b\\ a_j, \a_{[j]}\end{array}\!\!\right),
\end{equation*}
where $S(z)$ is the sum of residues of the integrand $\G(s)$ in \eqref{eq:G-defined} at the points where both $\Gamma(\a_1+s)$ and $\Gamma(\b_2+s)$ have poles (if no such points exist, then $S(z)=0$).  The function $G^{1,p}_{p,p}(w)$ on the right hand side can now be continued analytically to the sector $-\pi<\arg(w)<\pi$ by means of deformation of the integration contour  \cite[(11.35)]{Braaksma} or to the domain $|w|>1$ by the convergent sum of residues  \cite[(11.16), (11.36)]{Braaksma}. For $|w|<1$ this function can be expressed in terms of the generalized hypergeometric series times a power of $w$ \cite[(11.6)]{Braaksma}.  

It seems unnoticed, however, that the analytic continuations of  $\GI^{m,n}_{p,p}(z)$ and $\GO^{m,n}_{p,p}(z)$ in the ''balanced'' case $m+n=p$ and the real parameters take a particularly simple form expressed in terms of a single $G$ function. Unlike the expansion in hypergeometric functions, the formula obtained in this way is independent of simplicity of poles of the integrand. Deriving this formula is the goal of the following Section~2.

In the subsequent Section~3 we will explore the properties of this analytic continuation on the banks of the branch cut complementing our previous results for the  generalized hypergeometric functions \cite[Theorem~3.1]{KPITSF2017}.  In Section~4 we derive a transformation for the particular Meijer's $G$ function $G^{p,0}_{p,p}$ (the Meijer-N{\o}rlund function) whose $p-2$ lower parameters exceed the corresponding $p-2$ upper parameters by positive integers (implying that the integrand  $\G(s)$ is a ratio of products of two gamma functions times a polynomial in $s$).  This is motivated by the transformations of the generalized hypergeometric function with integral parameters differences described in their full generality by Miller and Paris in \cite{MP2013} and further studied by us in \cite{KPChapter2020,KPResults2019}.

One application of the particular case $G^{p,p}_{2p,2p}$ of $G^{m,n}_{2p,2p}$ with $m+n=2p$ is the following curious identity deduced by us in \cite[Theorem~4.1]{KPITSF2017}:
\begin{multline}\label{eq:sum-integral}
\sum\limits_{k=0}^{m}\left\{\frac{\Gamma(\a+k)\Gamma(\a+\alpha+\beta+m-k)}{\Gamma(\b+k)\Gamma(\b+\alpha+\beta+m-k)}
-\frac{\Gamma(\a+\alpha+k)\Gamma(\a+\beta+m-k)}{\Gamma(\b+\alpha+k)\Gamma(\b+\beta+m-k)}\right\}
\\
=\int\limits_{0}^{1}G^{p,p}_{2p,2p}\left(x\left|\begin{array}{l}\!\!1-\a-m,\b+\alpha+\beta\!\!\\\!\!\a+\alpha+\beta,1-\b-m\!\!\end{array}\right.\right)
\!\frac{(1-x^{m+1})(1-x^{\alpha})(1-x^{\beta})}{x^{\alpha+\beta+1}(1-x)}dx.
\end{multline} 
It holds for $m\in\N_0$ and $\a,\b\in\C^p$ satisfying $0<\Re(a_j)<\sum_{j=1}^{p}\Re(b_j)$ for all $j$ and positive $\alpha$,$\beta$.    The final Section~5 of this paper is devoted to the extension and application of this identity  
by replacing $2p$ with arbitrary integer $r>p$. Moreover, we show that this result actually follows from a simpler integral involving  $G^{p,p}_{r,r}$.  
Finally, we apply this integral for summing a series involving digamma function whose particular case gives a new formula for the power series coefficients of the product of two generalized hypergeometric functions with shifted parameters.  

\section{Analytic continuation}

Recall from \eqref{eq:notation} that the sine function  of a vector stands for the product of sines of its components.  We will need the following lemma, first discovered by N{\o}rlund \cite[(3.41)]{Norlund} who did not present a detailed proof. We give an independent derivation below.  
\begin{lemma}\label{lm:Hermite}
For any complex vectors $\a$, $\b$ of size $m$ the following identity holds
\begin{equation}\label{eq:hermite}
\sum\limits_{k=1}^{m}\frac{\sin(\pi(\b-a_k))}{\sin(\pi(\a_{[k]}-a_k))}
\frac{e^{\pm{i}\pi(z-a_k)}}{\sin(\pi(z-a_k))}=e^{\pm{i}\pi\psi_m}-\frac{\sin(\pi(z-\b))}{\sin(\pi(z-\a))},
\end{equation}
where $\psi_m=\sum_{j=1}^{m}(b_j-a_j)$.  
\end{lemma}
\textbf{Proof.} Indeed expanding $e^{\pm{i}\pi(z-a_k)}$ by Euler's formula we obtain two sums, the second of which reduces  to \cite[(4.2)]{Johnson},\cite[(3.14)]{KPSIGMA2016}
\begin{equation}\label{eq:ptolemy}
\sum\limits_{k=1}^{m}\frac{\sin(\pi(\b-a_k))}{\sin(\pi(\a_{[k]}-a_k))}
=\sin(\pi\psi_m),
\end{equation} 
while the first can be summed by the following generalization of Hermite's identity due to Johnson \cite[Theorem~5]{Johnson}
$$
\sum\limits_{k=1}^{m}\frac{\sin(\pi(\b-a_k))}{\sin(\pi(\a_{[k]}-a_k))}
\cot(\pi(z-a_k))=\cos(\pi\psi_m)-\frac{\sin(\pi(z-\b))}{\sin(\pi(z-\a))}.~~~~~\square
$$
Setting $z=b_j$, $j=1,\ldots,m$, in \eqref{eq:hermite}, the second summand on the right hand side vanishes and we obtain
\begin{equation}\label{eq:ptolemy-exp}
\sum\limits_{k=1}^{m}\frac{\sin(\pi(\b-a_k))}{\sin(\pi(\a_{[k]}-a_k))}
\frac{e^{\pm{i}\pi(b_j-a_k)}}{\sin(\pi(b_j-a_k))}=e^{\pm{i}\pi\psi_m}
\end{equation}
for any  $j=1,\ldots,m$.  The main result of this section is the following
\begin{theorem}\label{th:Gmnpp}
Suppose $m\ge1$, $m+n=p$ and $\a$, $\b$ are real vectors of size $p$. For any $z\in\C\setminus((-\infty,0]\cup[1,\infty))$ with $|z|>1$ analytic continuation of the function defined by the series \eqref{eq:sumresleft} is given by
\begin{multline}\label{eq:Gmnpp-continuation}
\GI^{m,n}_{p,p}\left(z\left|\begin{array}{l}\!\b\!\\\!\a\!\end{array}\right.\right)
=-\exp(\mp i\pi\psi_{m})\GI^{p,0}_{p,p}\!\left(\!\frac{1}{z}\left|\!\begin{array}{l}1-\a\\1-\b\end{array}\right.\!\!\!\right)
+\GI^{n,m}_{p,p}\left(\frac{1}{z}\left|\begin{array}{l}\!1-\a\!\\\!1-\b\!\end{array}\right.\right)
\\
=-\exp(\mp i\pi\psi_{m})\GO^{0,p}_{p,p}\!\left(\!z\left|\!\begin{array}{l}\b\\\a\end{array}\right.\!\!\!\right)
+\GO^{m,n}_{p,p}\left(z\left|\begin{array}{l}\!\b\!\\\!\a\!\end{array}\right.\right),
\end{multline}
where $\psi_m=\sum_{j=1}^{m}(b_{n+j}-a_j)$ and the ''$-$'' sign is chosen if $\Im(z)>0$ while the ''$+$'' sign is chosen if $\Im(z)<0$.
\end{theorem}
\textbf{Proof.} We will first assume that the components of $\b$ are distinct modulo integers.  Under this assumption, the formula for analytic continuation of the generalized hypergeometric function \cite[formula 4.2.4(22)]{KF}, \cite[formula (16.8.6)]{NIST} after simple renaming of parameters and application of the  reflection formula for the gamma function takes the form:
\begin{multline}\label{eq:pFp-1infinity-tr}
{_pF_{p-1}}\!\left(\begin{array}{l}1-\b+a_k\\1-\a_{[k]}+a_k\end{array}\vline\: z\right)
\\
=-\frac{\Gamma(1-\a_{[k]}+a_k)}{z\Gamma(1-\b+a_k)}\sum\limits_{j=1}^{p}\frac{\Gamma(b_{j}-\b_{[j]})}{\Gamma(b_{j}-\a)}
\frac{\pi(e^{\mp{i\pi}}z)^{b_j-a_k}}{\sin(\pi(b_j-a_k))}{_{p}F_{p-1}}\!\!\left(\!\!\begin{array}{l}1+\a-b_{j}\\1+\b_{[j]}-b_{j}\end{array}\vline\,\,\frac{1}{z}\!\right),~~~\Im(z)\gtrless0.
\end{multline}
Let $z=re^{i\phi}$, $-\pi<\phi<\pi$. If $\Im(z)>0$, then  $e^{-i\pi}z=re^{i(\phi-\pi)}$ and we get
\begin{equation}\label{eq:zcombined}
z^{a_k}(e^{-i\pi}z)^{b_{j}-a_{k}}=r^{a_k}r^{b_{j}-a_{k}}e^{i\phi{a_k}}e^{i\phi(b_j-a_k)}e^{-i\pi(b_j-a_k)}
=z^{b_j}(\cos(\pi(b_j-a_k))-i\sin(\pi(b_j-a_k))).
\end{equation}
Write $\a=(\a^{1},\a^{2})$, where $\a^{1}$ has $m$ component, $\a^{2}$ has $p-m=n$ components. In a similar fashion, $\b=(\b^{1},\b^{2})$, where $\b^{1}$ has $n$ components and $\b^{2}$ has $p-n=m$ components.  Then, according to \cite[(16.17.2)]{NIST} or \cite[8.2.2.3]{PBM3} for $|z|<1$, we have 
\begin{equation}\label{eq:Gmnpp-pFp-1}
\GI^{m,n}_{p,p}\left(z\left|\begin{array}{l}\!\b\!\\\!\a\!\end{array}\right.\right)
=\sum\limits_{k=1}^{m}z^{a_k}\frac{\Gamma(\a^{1}_{[k]}-a_k)\Gamma(1-\b^{1}+a_k)}{\Gamma(1-\a^{2}+a_k)\Gamma(\b^{2}-a_k)}
{}_{p}F_{p-1}\!\!\left(\!\!\begin{array}{l}1-\b+a_k\\1-\a_{[k]}+a_k\end{array}\!\vline\,z\!\right).
\end{equation}
The right hand side here is holomorphic in $\C\setminus((-\infty,0]\cup[1,\infty))$, and thus defines the analytic continuation of the left hand side for $|z|>1$. Then, in view of \eqref{eq:pFp-1infinity-tr}, \eqref{eq:zcombined} and \eqref{eq:ptolemy-exp}, we obtain:
\begin{multline}\label{eq:Gmnpp-expansion}
\GI^{m,n}_{p,p}\left(z\left|\begin{array}{l}\!\b\!\\\!\a\!\end{array}\right.\right)=-\frac{\pi}{z}\sum\limits_{k=1}^{m}\frac{\Gamma(\a^{1}_{[k]}-a_k)\Gamma(1-\b^{1}+a_k)\Gamma(1-\a_{[k]}+a_k)}
{\Gamma(1-\a^{2}+a_k)\Gamma(\b^{2}-a_k)\Gamma(1-\b+a_k)}
\\
\times\sum\limits_{j=1}^{p}\frac{\Gamma(b_{j}-\b_{[j]})}{\Gamma(b_{j}-\a)}
\frac{z^{a_k}(-z)^{b_j-a_k}}{\sin(\pi(b_j-a_k))}{_{p}F_{p-1}}\!\!\left(\!\!\begin{array}{l}1+\a-b_{j}\\1+\b_{[j]}-b_{j}\end{array}\vline\,\,\frac{1}{z}\!\right)
\\
=-\frac{\pi}{z}\sum\limits_{j=1}^{p}z^{b_j}\frac{\Gamma(b_{j}-\b_{[j]})}{\Gamma(b_{j}-\a)}
{_{p}F_{p-1}}\!\!\left(\!\!\begin{array}{l}1+\a-b_{j}\\1+\b_{[j]}-b_{j}\end{array}\vline\,\,\frac{1}{z}\!\right)
\sum\limits_{k=1}^{m}\frac{\Gamma(\a^{1}_{[k]}-a_k)\Gamma(1-\a^{1}_{[k]}+a_k)}
{\Gamma(\b^{2}-a_k)\Gamma(1-\b^{2}+a_k)}\frac{e^{-i\pi(b_j-a_k)}}{\sin(\pi(b_j-a_k))}
\\
=-\frac{1}{z}\sum\limits_{j=1}^{p}z^{b_j}\frac{\Gamma(b_{j}-\b_{[j]})}{\Gamma(b_{j}-\a)}
{_{p}F_{p-1}}\!\!\left(\!\!\begin{array}{l}1+\a-b_{j}\\1+\b_{[j]}-b_{j}\end{array}\vline\,\,\frac{1}{z}\!\right)
\left\{e^{-i\pi\psi_m}-\frac{\sin(\pi(b_j-\b^{2}))}{\sin(\pi(b_j-\a^{1}))}\right\}
\\
=-\frac{e^{-i\pi\psi_m}}{z}\sum\limits_{j=1}^{p}z^{b_j}\frac{\Gamma(b_{j}-\b_{[j]})}{\Gamma(b_{j}-\a)}
{_{p}F_{p-1}}\!\!\left(\!\!\begin{array}{l}1+\a-b_{j}\\1+\b_{[j]}-b_{j}\end{array}\vline\,\,\frac{1}{z}\!\right)
\\
+\frac{1}{z}\sum\limits_{j=1}^{n}z^{b_j}\frac{\Gamma(b_{j}-\b_{[j]})}{\Gamma(b_{j}-\a)}\frac{\sin(\pi(b_j-\b^{2}))}{\sin(\pi(b_j-\a^{1}))}
{_{p}F_{p-1}}\!\!\left(\!\!\begin{array}{l}1+\a-b_{j}\\1+\b_{[j]}-b_{j}\end{array}\vline\,\,\frac{1}{z}\!\right)
\\
=-e^{-i\pi\psi_m}\GI^{p,0}_{p,p}\!\left(\!\frac{1}{z}\left|\!\begin{array}{l}1-\a\\1-\b\end{array}\right.\!\!\!\right)
+\GI^{n,m}_{p,p}\left(\frac{1}{z}\left|\begin{array}{l}\!1-\a\!\\\!1-\b\!\end{array}\right.\right),
\end{multline}
where we used the fact that
$$
\frac{\sin(\pi(b_j-\b^{2}))}{\sin(\pi(b_j-\a^{1}))}=0
$$
for $j=n+1,\ldots,p$, Euler's reflection formula and the shifting property of $G$ function. This yields the $\Im(z)>0$ case of \eqref{eq:Gmnpp-continuation}.

If $\Im(z)<0$, then  $e^{i\pi}z=re^{i(\phi+\pi)}$. Thus for  the principal value of the power function we will have
$$
z^{a_k}(-z)^{b_{j}-a_{k}}=r^{a_k}r^{b_{j}-a_{k}}e^{i\phi{a_k}}e^{i\phi(b_j-a_k)}e^{i\pi(b_j-a_k)}
=z^{b_j}(\cos(\pi(b_j-a_k))+i\sin(\pi(b_j-a_k))).
$$
Then a similar calculation leads to the $\Im(z)<0$ case of \eqref{eq:Gmnpp-continuation}.
Finally, to remove the restriction that the components of $\b$ are distinct modulo integers, we note that both sides of \eqref{eq:Gmnpp-continuation} are continuous in $\b$.  
 $\hfill\square$
\smallskip

\textbf{Remark.}  The third line in \eqref{eq:Gmnpp-expansion} coincides with $l=1$ case of Meijer's formula \eqref{eq:MeijerContin}, while a similar calculation when $\Im(z)<0$ leads to $l=0$ case of Meijer's formula. For the sake of completeness we give full derivation in the proof without resorting to Meijer's result.  

As a corollary, we immediately obtain the analytic continuation for the external function $\GO^{m,n}_{p,p}$.
\begin{corollary}\label{th:Gmnpp}
Suppose $n\ge1$ and $m+n=p$. For any $z\in\C\setminus((-\infty,0]\cup[1,\infty))$ with $|z|<1$ analytic continuation of the function defined by \eqref{eq:sumresright} is given by
\begin{multline}\label{eq:Gmnpp-continuationIns}
\GO^{m,n}_{p,p}\!\left(z~\vline\begin{array}{l}\b\\\a\end{array}\!\right)
=-\exp(\pm i\pi\psi_{n})\GO^{0,p}_{p,p}\!\left(\!\frac{1}{z}\left|\!\begin{array}{l}1-\a\\1-\b\end{array}\right.\!\!\!\right)
+\GO^{n,m}_{p,p}\left(\frac{1}{z}\left|\begin{array}{l}\!1-\a\!\\\!1-\b\!\end{array}\right.\right)
\\
=-\exp(\pm i\pi\psi_{n})\GI^{p,0}_{p,p}\!\left(\!z\left|\!\begin{array}{l}\b\\\a\end{array}\right.\!\!\!\right)
+\GI^{m,n}_{p,p}\left(z\left|\begin{array}{l}\b\!\\\a\!\end{array}\right.\right),
\end{multline}
where $\psi_n=\sum_{j=1}^{n}(b_{j}-a_{m+j})$ and the ''$+$'' sign is chosen if $\Im(z)>0$ while the ''$-$'' sign is chosen if $\Im(z)<0$.
\end{corollary}
\textbf{Proof.} Indeed, for $|z|>1$ we can express the function of the left hand side of \eqref{eq:Gmnpp-continuationIns} in terms of $\GI^{n,m}_{p,p}(1/z)$ by \eqref{eq:reflection}.  We can then apply Theorem~\ref{th:Gmnpp} to get analytic continuation to $|z|<1$ which yields  \eqref{eq:Gmnpp-continuationIns}.  $\hfill\square$

\smallskip
Formula \eqref{eq:Gmnpp-continuation} simplifies significantly in the case of Meijer-N{\o}rlund function $\GI^{p,0}_{p,p}$. 
\begin{corollary}\label{th:Gp0pp}
For any $z\in\C\setminus((-\infty,0]\cup[1,\infty))$ with $|z|>1$ the analytic continuation is given by
\begin{equation}\label{eq:G-continuation}
\GI^{p,0}_{p,p}\!\left(\!z\left|\!\begin{array}{l}\b\\\a\end{array}\right.\!\!\!\right)
=-\exp(\mp i\pi\psi_{p})\GI^{p,0}_{p,p}\!\left(\!\frac{1}{z}\left|\!\begin{array}{l}1-\a\\1-\b\end{array}\right.\!\!\!\right)
=-\exp(\mp i\pi\psi_{p})\GO^{0,p}_{p,p}\!\left(\!z\left|\!\begin{array}{l}\b\\\a\end{array}\right.\!\!\!\right),
\end{equation}
where $\psi=\sum_{j=1}^{p}(b_j-a_j)$ and the ''$-$'' sign is chosen if $\Im(z)>0$ while the ''$+$'' sign is chosen if $\Im(z)<0$.
\end{corollary}
\textbf{Proof}. It is sufficient to note that for $|z|<1$ we have 
$$
\GI^{0,p}_{p,p}\!\left(\!z\left|\!\begin{array}{l}\b\\\a\end{array}\right.\!\!\!\right)
=G^{0,p}_{p,p}\!\left(\!z\left|\!\begin{array}{l}\b\\\a\end{array}\right.\!\!\!\right)
=0
$$
according to \cite[p.229, (5)]{Meijer} (see also \cite[(A.1)]{KPSIGMA2016}) as there are no poles inside the contour. Application of  \eqref{eq:Gmnpp-continuation} now yields the result.  $\hfill\square$

\textbf{Remark.} In view of the formula \cite[formula (8.4.49.22)]{PBM3} 
$$
G^{2,0}_{2,2}\!\left(\!z~\vline\begin{array}{l}b_1,b_2
\\a_1,a_2\end{array}\!\!\right)=\frac{z^{a_2}(1-z)_{+}^{b_1+b_2-a_1-a_2-1}}{\Gamma(b_1+b_2-a_1-a_2)}
{_2F_{1}}\!\left(\begin{array}{l}b_1-a_1,b_2-a_1
\\b_1+b_2-a_1-a_2\end{array}\!\!;1-z\right)
$$
the case $p=2$ of the  first identity in \eqref{eq:G-continuation} is easily seen to be equivalent to the classical Euler-Pfaff transformation 
$$
{_2F_{1}}\!\left(\begin{array}{c}a,b
\\c\end{array}\!\!;z\right)=(1-z)^{-a}{_2F_{1}}\!\left(\begin{array}{c}a,c-b
\\c\end{array}\!\!;\frac{z}{z-1}\right).
$$

\section{The principal branch of $\GI^{m,n}_{p,p}$ on the banks of the branch cut}

As we mentioned above, the function $\GI^{m,n}_{p,p}$ is single-valued in the entire complex plane cut along the rays $(-\infty,0]$ и $[1,+\infty)$. The points on the cuts will be understood as accessible boundary points of the domain $\C\setminus[(-\infty,0]\cup[1,+\infty)]$.   In our first  theorem we will give the values of $\GI^{m,n}_{p,p}(z)$ on the banks of the cut $(-1,0]$.
\begin{theorem}\label{th:Gp0pp}
For  $m+n=p$, $0<x<1,$ $x_{\pm}=xe^{\pm i\pi}$  and real $\a$, $\b$ we have 
\begin{equation}\label{eq:G-upperbank}
\Re\left( \GI^{m,n}_{p,p}\!\left(\!x_{+}\left|\!\begin{array}{l}\b\\\a\end{array}\right.\!\!\!\right)\right)=\Re \left(\GI^{m,n}_{p,p}\!\left(\!x_{-}\left|\!\begin{array}{l}\b\\\a\end{array}\right.\!\!\!\right)\right)
= -\pi \GI^{m,n}_{p+1,p+1}\!\left(\!x~\vline\begin{array}{l}\b,3/2
\\\a,3/2\end{array}\!\!\right)
\end{equation}
and
\begin{equation}\label{eq:G-lowerbank}
\Im\left( \GI^{m,n}_{p,p}\!\left(\!x_{+}\left|\!\begin{array}{l}\b\\\a\end{array}\right.\!\!\!\right)\right)=-\Im  \left(\GI^{m,n}_{p,p}\!\left(\!x_{-}\left|\!\begin{array}{l}\b\\\a\end{array}\right.\!\!\!\right)\right)
= \pi \GI^{m,n}_{p+1,p+1}\!\left(\!x~\vline\begin{array}{l}\b,1
\\\a,1\end{array}\!\!\right).
\end{equation}
\end{theorem}
\textbf{Proof.}  Application of \cite[8.2.2.3]{PBM3} ($=$ \eqref{eq:Gmnpp-pFp-1} for $m+n=p$) to each summand and once again to the resulting expression  yields 
\begin{multline}\label{eq:Gmnpp-expansion21}
\GI^{m,n}_{p,p}\left(x_{+}\left|\begin{array}{l}\!\b\!\\\!\a\!\end{array}\right.\right)+\GI^{m,n}_{p,p}\left(x_{-}\left|\begin{array}{l}\!\b\!\\\!\a\!\end{array}\right.\right)
\\
=2\sum\limits_{k=1}^{m}x^{a_k}\cos(\pi{a_k})\frac{\Gamma(\a^{1}_{[k]}-a_k)\Gamma(1-\b^{1}+a_k)}{\Gamma(1-\a^{2}+a_k)\Gamma(\b^{2}-a_k)}
{}_{p}F_{p-1}\!\!\left(\!\!\begin{array}{l}1-\b+a_k\\1-\a_{[k]}+a_k\end{array}\!\vline\,-x\!\right)
\\ 
=-2\pi\sum\limits_{k=1}^{m}x^{a_k}\frac{\Gamma(\a^{1}_{[k]}-a_k)\Gamma(1-\b^{1}+a_k)}{\Gamma(a_k-1/2)\Gamma(3/2-a_k)\Gamma(1-\a^{2}+a_k)\Gamma(\b^{2}-a_k)}
{}_{p}F_{p-1}\!\!\left(\!\!\begin{array}{l}1-\b+a_k\\1-\a_{[k]}+a_k\end{array}\!\vline\,-x\!\right)
\\
=-2\pi\GI^{m,n}_{p+1,p+1}\!\left(\!x~\vline\begin{array}{l}\b,3/2
\\\a,3/2\end{array}\!\!\right)
\end{multline}
leading immediately to \eqref{eq:G-upperbank}.   In a similar fashion,
\begin{multline}\label{eq:Gmnpp-expansion21}
\GI^{m,n}_{p,p}\left(x_{+}\left|\begin{array}{l}\!\b\!\\\!\a\!\end{array}\right.\right)-\GI^{m,n}_{p,p}\left(x_{-}\left|\begin{array}{l}\!\b\!\\\!\a\!\end{array}\right.\right)
\\
=2i\sum\limits_{k=1}^{m}x^{a_k}\sin(\pi{a_k})\frac{\Gamma(\a^{1}_{[k]}-a_k)\Gamma(1-\b^{1}+a_k)}{\Gamma(1-\a^{2}+a_k)\Gamma(\b^{2}-a_k)}
{}_{p}F_{p-1}\!\!\left(\!\!\begin{array}{l}1-\b+a_k\\1-\a_{[k]}+a_k\end{array}\!\vline\,-x\!\right)
\\ 
=2i\pi\sum\limits_{k=1}^{m}x^{a_k}\frac{\Gamma(\a^{1}_{[k]}-a_k)\Gamma(1-\b^{1}+a_k)}{\Gamma(a_k)\Gamma(1-a_k)\Gamma(1-\a^{2}+a_k)\Gamma(\b^{2}-a_k)}
{}_{p}F_{p-1}\!\!\left(\!\!\begin{array}{l}1-\b+a_k\\1-\a_{[k]}+a_k\end{array}\!\vline\,-x\!\right)
\\
=2i\pi\GI^{m,n}_{p+1,p+1}\!\left(\!x~\vline\begin{array}{l}\b,1
\\\a,1\end{array}\!\!\right)
\end{multline}
leading immediately to \eqref{eq:G-lowerbank}.  $\hfill\square$

\medskip

The values of the function $\GI^{m,n}_{p,p}(z)$ on the banks of the cut $(-\infty,-1)$ turn out to be more cumbersome:
\begin{theorem}\label{th:banksoutside}
For  $m+n=p$, $x_{\pm}=xe^{\pm i\pi}$  with $x>1,$ $\psi_m=\sum_{j=1}^{m}(b_{n+j}-a_j)$ and real $\a,\,\b$  we have 
\begin{equation}\label{eq:banksoutsideRe}
\Re\left(\GI^{m,n}_{p,p}\!\left(\!x_{\pm}\left|\!\begin{array}{l}\b\\\a\end{array}\right.\!\!\!\right)\right)=-\pi\GI^{p,0}_{p+1,p+1}\!\left(\!\frac{1}{x}\left|\!\begin{array}{l}1-\a,1/2-\psi_m\\1-\b,1/2-\psi_m\end{array}\right.\!\!\!\right)
-\pi\GI^{n,m}_{p+1,p+1}\left(\frac{1}{x}\left|\begin{array}{l}\!1-\a,3/2\!\\\!1-\b,3/2\!\end{array}\right.\right)
\end{equation}
and
\begin{equation}\label{eq:banksoutsideIm}
\Im\left(\GI^{m,n}_{p,p}\!\left(\!x_{\pm}\left|\!\begin{array}{l}\b\\\a\end{array}\right.\!\!\!\right)\right)
=\pm\pi\GI^{p,0}_{p+1,p+1}\!\left(\!\frac{1}{x}\left|\!\begin{array}{l}1-\a,1-\psi_m\\1-\b,1-\psi_m\end{array}\right.\!\!\!\right)\mp\pi \GI^{n,m}_{p+1,p+1}\left(\frac{1}{x}\left|\begin{array}{l}\!1-\a,1\!\\\!1-\b,1\!\end{array}\right.\right).
\end{equation}
\end{theorem}
\textbf{Proof.}  According to Theorem~\ref{th:Gmnpp} we have
\begin{equation}\label{eq:G-continuation13}
\GI^{m,n}_{p,p}\!\left(\!xe^{i\pi}\left|\!\begin{array}{l}\b\\\a\end{array}\right.\!\!\!\right)
= -\exp(- i\pi\psi_{m})\GI^{p,0}_{p,p}\!\left(\!\frac{e^{-i\pi}}{x}\left|\!\begin{array}{l}1-\a\\1-\b\end{array}\right.\!\!\!\right)
+\GI^{n,m}_{p,p}\left(\frac{e^{-i\pi}}{x}\left|\begin{array}{l}\!1-\a\!\\\!1-\b\!\end{array}\right.\right),
\end{equation}
while  Theorem~\ref{th:Gp0pp} yields:
$$ 
\GI^{n,m}_{p,p}\left(\frac{e^{-i\pi}}{x}\left|\begin{array}{l}\!1-\a\!\\\!1-\b\!\end{array}\right.\right)=
-\pi\GI^{n,m}_{p+1,p+1}\left(\frac{1}{x}\left|\begin{array}{l}\!1-\a,3/2\!\\\!1-\b,3/2\!\end{array}\right.\right)-\pi i\GI^{n,m}_{p+1,p+1}\left(\frac{1}{x}\left|\begin{array}{l}\!1-\a,1\!\\\!1-\b,1\!\end{array}\right.\right)
$$
and 
\begin{multline*}
 -\exp(-i\pi\psi_{m})\GI^{p,0}_{p,p}\!\left(\!\frac{e^{-i\pi}}{x}\left|\!\begin{array}{l}1-\a\\1-\b\end{array}\right.\!\!\!\right)
 \\
 =(-\cos(\pi\psi_m)+i\sin(\pi\psi_m))\left\{-\pi\GI^{p,0}_{p+1,p+1}\!\left(\!\frac{1}{x}\left|\!\begin{array}{l}1-\a,3/2\\1-\b,3/2\end{array}\right.\!\!\!\right)
 -\pi{i}\GI^{p,0}_{p+1,p+1}\!\left(\!\frac{1}{x}\left|\!\begin{array}{l}1-\a,1\\1-\b,1\end{array}\right.\!\!\!\right)\right\}
 \\
 =-\pi\GI^{p,0}_{p+1,p+1}\!\left(\!\frac{1}{x}\left|\!\begin{array}{l}1-\a,1/2-\psi_m\\1-\b,1/2-\psi_m\end{array}\right.\!\!\!\right)+\pi i\GI^{p,0}_{p+1,p+1}\!\left(\!\frac{1}{x}\left|\!\begin{array}{l}1-\a,1-\psi_m\\1-\b,1-\psi_m\end{array}\right.\!\!\!\right).
\end{multline*}
The last equality is derived by combining the definition \eqref{eq:G-defined} of $G$ function and standard trigonometric identities. So we arrive at \eqref{eq:banksoutsideRe} and \eqref{eq:banksoutsideIm}.$\hfill\square$

Note that the values of $\GI^{m,n}_{p,p}$ on the cut $[1,+\infty)$ can be derived directly from Theorem~\ref{th:Gmnpp} leading to the following result:
\begin{theorem}\label{th:Gp0pp1}
For  $m+n=p$, $x_{\pm}=xe^{\pm i0}$ with $x>1$ and real $\a$, $\b$ we have 
\begin{multline}\label{eq:G-Rebank1}
\Re\left( \GI^{m,n}_{p,p}\!\left(\!x_{+}\left|\!\begin{array}{l}\b\\\a\end{array}\right.\!\!\!\right)\right)=\Re \left(\GI^{m,n}_{p,p}\!\left(\!x_{-}\left|\!\begin{array}{l}\b\\\a\end{array}\right.\!\!\!\right)\right)
=-\cos( \pi\psi_{m})\GI^{p,0}_{p,p}\!\left(\!\frac{1}{x}\left|\!\begin{array}{l}1-\a\\1-\b\end{array}\right.\!\!\!\right)
\\
+\GI^{n,m}_{p,p}\left(\frac{1}{x}\left|\begin{array}{l}\!1-\a\!\\\!1-\b\!\end{array}\right.\right)
\end{multline}
and
\begin{equation}\label{eq:G-Imbank1}
\Im\left( \GI^{m,n}_{p,p}\!\left(\!x_{+}\left|\!\begin{array}{l}\b\\\a\end{array}\right.\!\!\!\right)\right)=-\Im\left(\GI^{m,n}_{p,p}\!\left(\!x_{-}\left|\!\begin{array}{l}\b\\\a\end{array}\right.\!\!\!\right)\right)
=\sin(\pi\psi_{m})\GI^{p,0}_{p,p}\!\left(\!\frac{1}{x}\left|\!\begin{array}{l}1-\a\\1-\b\end{array}\right.\!\!\!\right).
\end{equation}
\end{theorem}

\section{Transformations for $\GI_{p,p}^{p,0}$ with integral parameter differences}

In \cite[(3.1)]{KPITSF2017} we showed that for any $\a\in\R^p$, $\b\in\R^{p-1}$, $d\in \R$ and $0<x<1$ the following identity holds
  \begin{multline}\label{eq:GF}
\GI^{p,0}_{p,p}\left(x\,\vline\begin{array}{c}d,\b\\\a\end{array}\!\!\right)\!=\!
\frac{\Gamma(\a-d+1)x^{d-1}}{2\pi{i}\Gamma(\b-d+1)}
\left\{{_{p}F_{p-1}}\!\left(\!\!\begin{array}{l}\a-d+1\\\b-d+1\end{array}\vline\,\,\frac{1}{x}+i0\!\right)-
{_{p}F_{p-1}}\!\left(\!\!\begin{array}{l}\a-d+1\\\b-d+1\end{array}\vline\,\,\frac{1}{x}-i0\!\right)\right\}.
\end{multline}
 Relation \eqref{eq:GF} enables us to transfer some of the properties of the hypergeometric functions to the Meijer-N{\o}rlund function $\GI^{p,0}_{p,p}$.  In this section we will illustrate this idea by deriving an analogue of the Miller-Paris transformations valid for the generalized hypergeometric functions with integral parameters differences. 
 These transformations are described below.  Let $\m=(m_1,\ldots,m_r)\in\N^r$ ($r$-tuple of natural numbers), $m=m_1+m_2+\ldots+m_r$, $\f=(f_1,\ldots,f_r)\in\C^r$. The most general Miller--Paris transformations are given in \cite[Theorems~3~and~4]{MP2013}. The following identity  generalizes the first Euler (or Euler--Pfaff) transformation for the Gauss hypergeometric function:
\begin{equation}\label{eq:KRPTh1-1}
{}_{r+2}F_{r+1}\left.\!\!\left(\begin{matrix}\aaaa, \bbbb, \f+\m\\\cccc,\f\end{matrix}\right\vert x\right)
=(1-x)^{-\aaaa}{}_{m+2}F_{m+1}\left.\!\!\left(\begin{matrix}\aaaa, \cccc-\bbbb-m, \zetta+1\\\cccc, \zetta\end{matrix}\right\vert\frac{x}{x-1}\right).
\end{equation}
It is true for $\bbbb\ne{f_j}$, $j=1,\ldots,r$,  $(\cccc-\bbbb-m)_{m}\ne0$ and $x\in\C\setminus[1,\infty)$, where as usual the generalized hypergeometric function is understood as the analytic continuation for the values of the argument making the series definition diverge.  Here  $\zetta=\zetta(\bbbb,\cccc,\f)=(\zeta_1,\ldots,\zeta_m)$ are the roots of the characteristic polynomial
\begin{equation}\label{eq:Qm}
Q_m(t)=Q_m(\bbbb,\cccc,\f|t)=\frac{1}{(\cccc-\bbbb-m)_{m}}\sum\limits_{k=0}^{m}(\bbbb)_k(t)_{k}(\cccc-\bbbb-m-t)_{m-k}
\frac{(-1)^k}{k!}{}_{r+1}F_{r}\!\left(\begin{matrix}-k,\f+\m\\\f\end{matrix}\right).
\end{equation}
The other transformation generalizes the second Euler transformation for the Gauss
hypergeometric function as follows:
\begin{equation}\label{eq:KRPTh1-2}
{}_{r+2}F_{r+1}\left.\!\!\left(\begin{matrix}\aaaa, \bbbb, \f+\m\\\cccc,\f\end{matrix}\right\vert x\right)
=(1-x)^{\cccc-\aaaa-\bbbb-m}{}_{m+2}F_{m+1}\left.\!\!\left(\begin{matrix}\cccc-\aaaa-m, \cccc-\bbbb-m, \etta+1\\\cccc, \etta\end{matrix}\right\vert x\right)
\end{equation}
and holds true for $(\cccc-\aaaa-m)_{m}\ne0$, $(\cccc-\bbbb-m)_{m}\ne0$, $(1+\aaaa+\bbbb-\cccc)_m\ne0$ and $x\in\C\setminus[1,\infty)$.
The vector $\etta=\etta(\aaaa,\bbbb,\cccc,\f)=(\eta_1,\ldots,\eta_m)$ is formed by the roots of the second characteristic polynomial $\hat{Q}_m(t)=\hat{Q}_m(a,b,c,\f|t)$ given by
\begin{equation}\label{eq:hatP}
\hat{Q}_m(t)=\sum\limits_{k=0}^m\frac{(-1)^k(a)_k(-b-m)_k(t)_k(c-a-m-t)_{m-k}}{(c-a-m)_m(c-b-m)_kk!}
{}_{r+2}F_{r+1}\!\left(\begin{matrix}-k,b,\f+\m\\b+m-k+1,\f\end{matrix}\right).
\end{equation}
Note that this form of the characteristic polynomial established by us in \cite[Theorem~2]{KPChapter2020} is simpler than the original form given in \cite{MP2013}. 
Transformations \eqref{eq:KRPTh1-1}, \eqref{eq:KRPTh1-2} fail when $(\cccc-\bbbb-m)_{m}=0$, that is when $\cccc-\bbbb\in\{1,\ldots,m\}$.   Their  degenerate forms under this condition were established by us in \cite{KPChapter2020} and \cite{KPResults2019}.
These transformations lead to the following result for the Meijer-N{\o}rlund function:
\begin{theorem}\label{th:GParis2}
Suppose  $a,\, b,\,c,\,d\in\C$,  $\f\in\C^{p-2}$, $\m\in\N^{p-2}$, $m=m_1+\cdots+m_{p-2}$. Then the following identities hold
\begin{multline}
\GI^{p,0}_{p,p}\left(z\,\vline\begin{array}{c}c,d,\f\\a,b,\f+\m\end{array}\!\!\right)=
z^{d-\alpha-1}(1-z)^\alpha\frac{\Gamma(\a-d+1)\Gamma(\hat{\b\!\!}\,\,)}{\Gamma(\hat{\a})\Gamma(\b-d+1)}
\\
\times\left\{\cos(\pi\alpha)\GI^{m+2,0}_{m+2,m+2}\left(z\,\vline\begin{array}{c}1,\hat{\b\!\!}\,\\\hat{\a}\end{array}\!\!\right)
+\sin(\pi\alpha)
\GI^{m+2,1}_{m+3,m+3}\left(z\,\vline\begin{array}{c}1,3/2,\hat{\b\!\!}\,\\\hat{\a},3/2\end{array}\!\!\right)\right\}
\end{multline}
and
\begin{multline}
\GI^{p,0}_{p,p}\left(z\,\vline\begin{array}{c}c,d,\f\\a,b,\f+\m\end{array}\!\!\right)=
z^{d-\beta-1}(1-z)^\beta\frac{\Gamma(\a-d+1)\Gamma(\bar{\b\!\!}\,\,)}{\Gamma(\bar{\a})\Gamma(\b-d+1)}
\\
\times\left\{-\cos(\pi\beta)\GI^{m+2,0}_{m+2,m+2}\left(1-z\,\vline\begin{array}{c}1,\bar{\b\!\!}\,\\\bar{\a}\end{array}\!\!\right)
+\sin(\pi\beta)
\GI^{m+2,1}_{m+3,m+3}\left(1-z\,\vline\begin{array}{c}1,3/2,\bar{\b\!\!}\,\\\bar{\a},3/2\end{array}\!\!\right)\right\},
\end{multline}
where  $\a=(a,b,\f+\m)$, $\b=(c,\f)$, $\alpha=c+d-a-b-m-1$, $\hat{\a}=(c-a-m,c-b-m,{\lambdda}+1)$, $\hat{\b\!\!}\,=(c-d+1,{\lambdda})$ and $\lambdda=(\lambda_1,\ldots,\lambda_m)$ are the roots of the second characteristic polynomial
$$
\hat{Q}_m(\aaaa-d+1,\bbbb-d+1,\cccc-d+1,\f-d+1|t),
$$ 
defined \eqref{eq:hatP}, $\beta=d-a-1$, $\bar{\a}=(a-d+1,c-b-m,{\ggamma}+1)$, $\bar{\b\!\!}\,=(c-d+1,{\ggamma})$ and
$\ggamma$ are the roots of the first characteristic polynomial
$$
{Q}_m(\bbbb-d+1,\cccc-d+1,\f-d+1|t)
$$
defined in \eqref{eq:Qm}.
\end{theorem}

\textbf{Proof.} Suppose  $a$, $b$, $c$, $d$ and $\f$ are real, $r=p-2$.  In \cite[Theorem~3.1]{KPITSF2017} we have shown that for  $x>1$
\begin{equation}\label{eq:upper}
{_{p+1}F_p}\!\left(\!\!\begin{array}{l}\a\\\b\end{array}\vline\,\,x+i0\!\right)
=-\frac{\pi\Gamma(\b)}{\sqrt{x}\Gamma(\a)}
G^{p+1,1}_{p+2,p+2}\left(\frac{1}{x}\,\vline\begin{array}{c}1/2,1,\b-1/2\\\a-1/2,1\end{array}\!\!\right)+
\pi{i}\frac{\Gamma(\b)}{\Gamma(\a)}G^{p+1,0}_{p+1,p+1}\left(\frac{1}{x}\,\vline\begin{array}{c}1,\b\\\a\end{array}\!\!\right)
\end{equation}
and
\begin{equation}\label{eq:low}
{_{p+1}F_p}\!\left(\!\!\begin{array}{l}\a\\\b\end{array}\vline\,\,x-i0\!\right)
=-\frac{\pi\Gamma(\b)}{\sqrt{x}\Gamma(\a)}
G^{p+1,1}_{p+2,p+2}\left(\frac{1}{x}\,\vline\begin{array}{c}1/2,1,\b-1/2\\\a-1/2,1\end{array}\!\!\right)-
\pi{i}\frac{\Gamma(\b)}{\Gamma(\a)}G^{p+1,0}_{p+1,p+1}\left(\frac{1}{x}\,\vline\begin{array}{c}1,\b\\\a\end{array}\!\!\right).
\end{equation}
To make notation more concise we will skip the dimension indicies from the notation of the hypergeometric function. Applying \eqref{eq:GF}, \eqref{eq:KRPTh1-2}, \eqref{eq:upper} and \eqref{eq:low}, for $0<x<1$ we will get 
\begin{multline}
G^{p,0}_{p,p}\left(x\,\vline\begin{array}{c}d,\b\\\a\end{array}\!\!\right)=
\frac{\Gamma(\a-d+1)}{2\pi i\Gamma(\b-d+1)}x^{d-1}\left\{(1-(1/x+i0))^\alpha{F}\!\left(\!\!\begin{array}{l}\hat{\a}\\\hat{\b\!\!\!}\,\,\end{array}\vline\,\,1/x+i0\!\right)
-\right.\\
\left.((1-(1/x-i0))^\alpha{F}\!\left(\!\!\begin{array}{l}\hat{\a}\\\hat{\b\!\!\!}\,\,\end{array}\vline\,\,1/x-i0\!\right)\right\}=(1-x)^\alpha x^{d-\alpha-1}\frac{\Gamma(\a-d+1)\Gamma(\hat{\b\!\!\!}\,\,\,)}{\Gamma(\hat{\a})\Gamma(\b-d+1)}\times\\\left\{(\cos\pi\alpha) G^{m+2,0}_{m+2,m+2}\left(x\,\vline\begin{array}{c}1,\hat{\b\!\!\!}\,\,\\\hat{\a}\end{array}\!\!\right)
+(\sin\pi\alpha)
G^{m+2,1}_{m+3,m+3}\left(x\,\vline\begin{array}{c}1,3/2,\hat{\b\!\!\!}\,\,\\\hat{\a},3/2\end{array}\!\!\right)\right\},
\end{multline}
where  $\alpha$, $\a$, $\b$,  $\hat{\a}$, $\hat{\b\!\!\!}\,\,$ and  $\lambdda$ are as defined in the theorem. 

If instead of \eqref{eq:KRPTh1-2} we apply transformation \eqref{eq:KRPTh1-1} to the identity \eqref{eq:GF}, we will get
\begin{multline}
G^{p,0}_{p,p}\left(x\,\vline\begin{array}{c}d,\b\\\a\end{array}\!\!\right)=
\frac{\Gamma(\a-d+1)}{2\pi{i}\Gamma(\b-d+1)}x^{d-1}\left\{(1-(1/x+i0))^\beta
{F}\!\left(\!\!\begin{array}{l}\bar{\a}\\\bar{\b\!\!}\,\end{array}\vline\,\,\frac{1/x+i0}{1/x-1+i0}\!\right)
\right.
\\
-\left.((1-(1/x-i0))^\beta{F}\!\left(\!\!\begin{array}{l}\bar{\a}\\\bar{\b\!\!}\,\end{array}\vline\,\,\frac{1/x-i0}{1/x-1-i0}\!\right)\right\}
=(1-x)^\beta z^{d-\beta-1}\frac{\Gamma(\a-d+1)\Gamma(\bar{\b\!\!}\,\,)}{\Gamma(\bar{\a})\Gamma(\b-d+1)}
\\
\times\left\{-(\cos\pi\beta) G^{m+2,0}_{m+2,m+2}\left(1-x\,\vline\begin{array}{c}1,\bar{\b\!\!}\,\\\bar{\a}\end{array}\!\!\right)
+(\sin\pi\beta)G^{m+2,1}_{m+3,m+3}\left(1-x\,\vline\begin{array}{c}1,3/2,\bar{\b\!\!}\,\\\bar{\a},3/2\end{array}\!\!\right)\right\}.
\end{multline} 
The condition on $x$ and reality of $d$, $\a$, $\b$, $\f$ can be removed according to the principle of analytic continuation. 
 $\hfill\square$

\section{A curious integral of $G$ function and its consequences}

We will write $\min(\a)$ for the minimal component of a real vector $\a$.
\begin{theorem}\label{thm:1}
Suppose $\a\in\C^{p}$, $\b\in\C^{q}$, $r=p+q$, $2p\ge{r},$  $\min(\Re(\a+\mu),\Re(\a+\nu))>0$ and if $p-q$ are even, we additionally require  that  \begin{equation}\label{eq:convergenceat1}
\Re\Big(\sum\nolimits_{k=1}^q b_k-\sum\nolimits_{j=1}^p a_j+(p-q)(1-\mu-\nu)/2\Big)>0. 
\end{equation}
Then 
\begin{equation}\label{eq:Gintegral}
\frac{\Gamma(\a+\mu)\Gamma(\a+\nu)}{\Gamma(\b+\mu)\Gamma(\b+\nu)}=\int\limits_{0}^{1}(x^{-\mu-1}+x^{-\nu-1})
G^{p,p}_{r,r}\left(x\left|\begin{array}{l}\!\!1-\a,\b+\mu+\nu\!\!\\\!\!\a+\mu+\nu,1-\b\!\end{array}\right.\right)dx.
\end{equation}
\end{theorem}
\textbf{Proof}. Suppose $\cb$, $\d$ are complex vectors of size $r=p+q$.  We will write $\d=(\d_{1},\d_{2})$, where $\d_{1}$ stands for the first $p$ components of $\d$ and $\d_{2}$ for the ultimate $q$ components.  Similarly, $\cb=(\cb_{1},\cb_{2})$ with $p$ components in $\cb_{1}$ and $q$ in $\cb_{2}$.   In this case, combining \cite[8.2.2.3]{PBM3} with \cite[8.2.2.4]{PBM3} and we can write for $x>0$:
\begin{multline}\label{eq:Gppexpansion}
G^{p,p}_{r,r}\left(x\left|\begin{array}{l}\!\!\cb_{1},\cb_{2}\!\!\\\!\!\d_{1},\d_{2}\!\!\end{array}\right.\right)
=H(1-x)\sum\limits_{j=1}^{p}A_jx^{d_j}{}_{r}F_{r-1}\!\!\left(\!\!\begin{array}{l}1-\cb+d_j\\1-\d_{[j]}+d_j\end{array}\!\vline\,(-1)^{p-q}x\!\right)
\\
+H(x-1)\sum\limits_{j=1}^{p}B_jx^{c_j-1}{}_{r}F_{r-1}\!\!\left(\!\!\begin{array}{l}1+\d-c_j\\1+\cb_{[j]}-c_j\end{array}\!\vline\,\frac{(-1)^{p-q}}{x}\!\right),
\end{multline}
where
$$
A_j=\frac{\Gamma(\d_{1{[j]}}-d_j)\Gamma(1-\cb_{1}+d_j)}{\Gamma(1-\d_{2}+d_j)\Gamma(\cb_{2}-d_j)},
~~~B_j=\frac{\Gamma(c_j-\cb_{1[j]})\Gamma(1+\d_{1}-c_j)}{\Gamma(1+\cb_{2}-c_j)\Gamma(c_j-\d_{2})}
$$
and $H(t)=\partial_t\max\{t,0\}$ is the Heaviside function. 
Next, according to \cite[2.24.2.1]{PBM3} (see also \cite[Th.~2.2 and 3.3]{KilSaig}):
\begin{equation*}
\frac{\Gamma(1-\cb_{1}-s)\Gamma(\d_{1}+s)}{\Gamma(\cb_{2}+s)\Gamma(1-\d_{2}-s)}=\int\limits_{0}^{\infty}x^{s-1}
G^{p,p}_{r,r}\left(x\left|\begin{array}{l}\!\!\cb\!\!\\\!\!\d\!\!\end{array}\right.\right)dx.
\end{equation*}

Suppose now that $\d=-\cb+u$, $u\in\C$ i.e. coordinate-wise $d_j=-c_j+u$, $j=1,\ldots,r$. Then it is straightforward to see that $A_j=B_j$ and the hypergeometric functions in the first and second term in \eqref{eq:Gppexpansion} coincide, leading to 
\begin{equation*}
\frac{\Gamma(1-\cb_{1}-s)\Gamma(u-\cb_{1}+s)}{\Gamma(\cb_{2}+s)\Gamma(1-u+\cb_{2}-s)}=
\int\limits_{0}^{1}(x^{s-1}+x^{-s-u})
G^{p,p}_{r,r}\left(x\left|\begin{array}{l}\!\!\a\!\!\\\!\!u-\a\!\!\end{array}\right.\right)dx.
\end{equation*}
Substituting $\cb_{1}=1-\a$, $\cb_{2}=\b+\mu+\nu$, $s=-\mu$, $u=\nu+\mu+1$ we obtain \eqref{eq:Gintegral}.  Convergence conditions are derived as follows: write $G$ function as a sum of hypergeometric functions according to  \eqref{eq:Gppexpansion}. Conditions $\Re(\a+\mu)>0$ and $\Re(\a+\nu) >0$ are needed to ensure convergence in the neighborhood of the origin. If $p-q$ is even the point $x=(-1)^{p-q}=1$ is singular. Using the expansion \cite[(31)]{CKP} of the hypergeometric function around this point (see also \cite[Theorem 2]{Buhring}), we see that the additional condition  \eqref{eq:convergenceat1} must be imposed.$\hfill\square.$

In view of the standard evaluation \cite[(4)]{CKP}, formula \eqref{eq:Gintegral} for $p=q$ can be rewritten as follows:
\begin{corollary}
Suppose $\a,\b\in\C^{p}$, $\sum\nolimits_{j=1}^p\Re(b_j-a_j)>0$ and $\min(\Re(\a+\mu),\Re(\a+\nu))>0$. Then 
\begin{multline*}
\int_{0}^{1}x^{\mu-1}G^{p,0}_{p,p}\left(x\left|\begin{array}{l}\!\!\b\!\!\\\!\!\a\end{array}\right.\right)dx\int_{0}^{1}x^{\nu-1}G^{p,0}_{p,p}\left(x\left|\begin{array}{l}\!\!\b\!\!\\\!\!\a\end{array}\right.\right)dx
\\
=\int_{0}^{1}(x^{-\mu-1}+x^{-\nu-1})
G^{p,p}_{2p,2p}\left(x\left|\begin{array}{l}\!\!1-\a,\b+\mu+\nu\!\!\\\!\!\a+\mu+\nu,1-\b\!\end{array}\right.\right)dx.
\end{multline*}
\end{corollary}

\begin{corollary}
Formula \eqref{eq:sum-integral} remains valid if $\a$, $\b$ have different sizes,
say if $\a\in\C^p$, $\b\in\C^q$.  In other words, writing $r=p+q$, $p\ge{q}$, we will have
\begin{multline}\label{eq:sum-integral-modified}
\sum\limits_{k=0}^{m}\left\{\frac{\Gamma(\a+k)\Gamma(\a+\alpha+\beta+m-k)}{\Gamma(\b+k)\Gamma(\b+\alpha+\beta+m-k)}
-\frac{\Gamma(\a+\alpha+k)\Gamma(\a+\beta+m-k)}{\Gamma(\b+\alpha+k)\Gamma(\b+\beta+m-k)}\right\}
\\
=\int_{0}^{1}G^{p,p}_{r,r}\left(x\left|\begin{array}{l}\!\!1-\a-m,\b+\alpha+\beta\!\!\\\!\!\a+\alpha+\beta,1-\b-m\!\!\end{array}\right.\right)
\!\frac{(1-x^{m+1})(1-x^{\alpha})(1-x^{\beta})}{x^{\alpha+\beta+1}(1-x)}dx.
\end{multline} 
\end{corollary} 
Finally, we apply Theorem~\ref{thm:1} to give a presumably new summation formula related to the power series coefficients of a product of two generalized hypergeometric series with shifted parameters. 
\begin{corollary}\label{th:linearization}
Suppose $p\geq q$, $\a\in\C^{p}$ and $\b\in\C^{q}$ and
\begin{equation}\label{eq:cond1}
\sum\nolimits_{j=1}^{q}b_j-\sum\nolimits_{i=1}^{p}a_j>\frac{1}{2}\left((p-q)(\lambda+m-1)+1\right).    
\end{equation}
Then we have the summation  formula
\begin{multline}\label{eq:psisum}
\sum\limits_{j=1}^{p}\frac{\pi^{p-q}\sin(\pi(\b-a_j))}{\sin(\pi(\a_{[j]}-a_j))}\sum\limits_{n=0}^\infty\frac{\Gamma(1-\b+a_j+n)\Gamma(\a+a_j+\lambda+m+n)(-1)^{(p-q)n}}{\Gamma(1-\a_{[j]}+a_j+n)\Gamma(\b+a_j+\lambda+m+n)n!}d_{j,n}(m)
\\
=\sum_{k=0}^{m}\frac{\Gamma(\a+k)\Gamma(\a+\lambda+m-k)}{\Gamma(\b+k)\Gamma(\b+\lambda+m-k)}=:a_m
\end{multline}
for $m=0,1,\ldots$, where \emph{(}$\psi(z)=\Gamma'(z)/\Gamma(z)$ denotes the digamma function\emph{)} 
$$
d_{j,n}(m)=\psi(a_j+\lambda+m+n+1)+\psi(a_j+m+n+1)-\psi(a_j+\lambda+n)-\psi(a_j+n).
$$
In particular, for $p=q$ we get the expansion
$$
{_{p+1}F_{q}}\left(\left.\!\!\begin{array}{c}\a,1\\\b\end{array}\right|z\!\right)
{_{p+1}F_{q}}\left(\left.\!\!\begin{array}{c}\a+\lambda,1\\\b+\lambda\end{array}\right|z\!\right)=\frac{\Gamma(\b)\Gamma(\b+\lambda)}{\Gamma(\a)\Gamma(\a+\lambda)}\sum\limits_{m=0}^{\infty} a_m z^m
$$
with $a_m$ as given above.
\end{corollary}
\textbf{Proof}. By using Theorem~\ref{thm:1} we obtain
\begin{multline}\label{eq:am}
a_m=\sum\limits_{k=0}^{m}\frac{\Gamma(\a+k)\Gamma(\a+\lambda+m-k)}{\Gamma(\b+k)\Gamma(\b+\lambda+m-k)}
\\
=\int\limits_{0}^{1}G^{p,p}_{p+q,p+q}\left(x\left|\begin{array}{l}\!\!-\a,\b+\lambda+m-1\!\!\\\!\!\a+\lambda+m-1,-\b\!\end{array}\right.\right)\sum\limits_{k=0}^{m}\left(x^{-k}+x^{k-m-\lambda}\right)dx
\\
=\int\limits_{0}^{1}\frac{(1-x^{m+1})}{x^m(1-x)}(1+x^{-\lambda})G^{p,p}_{p+q,p+q}\left(x\left|\begin{array}{l}\!\!-\a,\b+\lambda+m-1\!\!\\\!\!\a+\lambda+m-1,-\b\!\end{array}\right.\right)dx
\\
=\int\limits_{0}^{1}\frac{(1-x^{m+1})}{(1-x)}(1+x^{-\lambda})G^{p,p}_{p+q,p+q}\left(x\left|\begin{array}{l}\!\!-\a-m,\b+\lambda-1\!\!\\\!\!\a+\lambda-1,-\b-m\!\end{array}\right.\right)dx.
\end{multline}
We can write for $0<x<1$, $r=p+q$:
\begin{equation*}\label{eq:Gppexpansion2}
G^{p,p}_{r,r}\left(x\left|\begin{array}{l}\!\!-\a-m,\b+\lambda-1\!\!\\\!\!\a+\lambda-1,-\b-m\!\end{array}\right.\right)
\!=\!\sum\limits_{j=1}^{p}A_{j,m}x^{a_j+\lambda-1}{}_{r}F_{r-1}\!\!\left(\!\!\begin{array}{l}\a+a_j+\lambda+m,1-\b+a_j\\1-\a_{[j]}+a_j,\b+a_j+\lambda+m\end{array}\!\vline\,(-1)^{p-q} x\!\right),
\end{equation*}
where
$$
A_{j,m}=\frac{\Gamma(\a_{[j]}-a_j)\Gamma(\a+a_j+\lambda+m)}{\Gamma(\b-a_j)\Gamma(\b+a_j+\lambda+m)}=\frac{\pi^{p-q}\sin(\pi(\b-a_j))\Gamma(1-\b+a_j)\Gamma(\a+a_j+\lambda+m)}{\sin(\pi(\a_{[j]}-a_j))\Gamma(1-\a_{[j]}+a_j)\Gamma(\b+a_j+\lambda+m)}.
$$
It then follows from \eqref{eq:am} that
\begin{multline*}\label{eq:am1}
a_m=\sum\limits_{j=1}^{p}\frac{\pi^{p-q}\sin(\pi(\b-a_j))}{\sin(\pi(\a_{[j]}-a_j))}\sum\limits_{n=0}^\infty\frac{\Gamma(1-\b+a_j+n)\Gamma(\a+a_j+\lambda+m+n)(-1)^{(p-q)n}}{\Gamma(1-\a_{[j]}+a_j+n)\Gamma(\b+a_j+\lambda+m+n)n!}d_{j,n}(m),
\end{multline*}
where 
$$
d_{j,n}(m)=\int\limits_{0}^{1}\frac{x^{a_j+\lambda+n-1}(1-x^{m+1})dx}{1-x}+\int\limits_{0}^{1}\frac{x^{a_j+n-1}(1-x^{m+1})dx}{1-x}.
$$
By using the integral evaluation
$$
\int_0^1\frac{x^a(1-x^{m+1})dx}{1-x}=\psi(a+m+2)-\psi(a+1),
$$ 
we get
$$
d_{j,n}(m)=\psi(a_j+\lambda+m+n+1)+\psi(a_j+m+n+1)-\psi(a_j+\lambda+n)-\psi(a_j+n).
$$
It is easy to verify using Stirling's formula and the asymptotic relation 
$$
\psi(x)=\log x+O(1/x), \ x\to\infty,
$$
that condition \eqref{eq:cond1} is sufficient for absolute convergence of the series on the left hand side of \eqref{eq:psisum}. 
$\hfill\square$

\bigskip
\bigskip

\textbf{Acknowledgements.} This work has been supported by the Ministry of Science and
Higher Education of the Russian Federation (agreement No. 075-02-2021-1395). The second named author was also supported by RFBR (project 20-01-00018).

\end{document}